\documentclass{birkjour}
\usepackage{amsmath, amsthm, amsfonts, mathdots}
\usepackage{amssymb, dsfont}
\usepackage{amsbsy}
\usepackage{mathrsfs}

\input amssym.def
\input amssym
\usepackage{enumerate}
\usepackage{hyperref}
\hypersetup{
    pdftitle={Closed operator functional calculus in Banach modules and applications},    
    pdfauthor={Anatoly Baskakov}{Ilya Krishtal},     
    pdfsubject={Closed operator functional calculus},   
    pdfcreator={Ilya Krishtal},   
    pdfproducer={Ilya Krishtal}, 
    pdfkeywords={Banach modules}{Beurling spectrum}{Wiener's lemma}, 
     colorlinks,%
    citecolor=blue,%
    filecolor=black,%
    linkcolor=magenta,%
    urlcolor=black
}

\chardef\bslash=`\\ 

\makeatletter
\def\verbatim{\interlinepenalty\@M \@verbatim
   \leftskip\@totalleftmargin\advance\leftskip2pc
   \frenchspacing\@vobeyspaces \@xverbatim}
\makeatother
\newtheorem{thm}{Theorem}[section]
\newtheorem{cor}[thm]{Corollary}
\newtheorem{lem}[thm]{Lemma}
\newtheorem{prop}[thm]{Proposition}

\newtheorem{con}{Conclusion}[section]
\theoremstyle{definition}
\newtheorem{defn}{Definition}[section]

\theoremstyle{remark}
\newtheorem{rem}{Remark}[section]
\newtheorem{exmp}{Example}[section]

\numberwithin{equation}{section}

\newcommand{\begeq}{\begin {equation}}
\newcommand{\eq}{\end{equation}}
\newcommand{\bs}{\begin {split}}
\newcommand{\es}{\end{split}}
\newcommand{\bp}{\begin {prop}}
\newcommand{\ep}{\end {prop}}
\newcommand{\bt}{\begin {thm}}
\newcommand{\et}{\end {thm}}
\newcommand{\bcon}{\begin {con}}
\newcommand{\econ}{\end {con}}
\newcommand{\bc}{\begin {cor}}
\newcommand{\ec}{\end {cor}}
\newcommand{\bl}{\begin {lem}}
\newcommand{\el}{\end {lem}}
\newcommand{\bpf}{\begin {proof}}
\newcommand{\epf}{\end {proof}}
\newcommand{\bi}{\begin {itemize}}
\newcommand{\ei}{\end {itemize}}
\newcommand{\ben}{\begin {enumerate}}
\newcommand{\een}{\end {enumerate}}
\newcommand{\brem}{\begin {rem}}
\newcommand{\erem}{\end {rem}}

\newcommand{\bd}{\begin {defn}}
\newcommand{\ed}{\end {defn}}

\newcommand{\bex}{\begin {exmp}}
\newcommand{\eex}{\end {exmp}}
\newcommand{\bai}{b.a.i.~}
\newcommand{\cfbai}{cf-b.a.i.~}
\newcommand{\id}{\mathrm{id}}
\newcommand{\eye}{\mathds 1}
\newcommand{\AP}{\mathcal{AP}}

\newcommand{\norm}[1]{\left\|{#1} \right\|}

\newcommand{\eps}{\epsilon}

\newcommand{\A}{\mathcal A}

\newcommand{\F}{\mathcal{F}}

\newcommand{\Y}{\mathcal{Y}}
\newcommand{\TTT}{{\T\kern-.44em \T}}
\newcommand{\tTTT}{\widetilde{\T\kern-.44em \T}}
\newcommand{\chT}{\check{\T}}

\newcommand{\ZZ}{{\mathbb Z}}
\newcommand{\TT}{{\mathbb T}}
\newcommand{\RR}{{\mathbb R}}
\newcommand{\CC}{{\mathbb C}}
\newcommand{\NN}{{\mathbb N}}

\newcommand{\X}{\mathcal{X}}

\newcommand{\s}{\sigma}
\renewcommand{\l}{\lambda}
\renewcommand{\L}{\Lambda}
\renewcommand{\a}{\alpha}

\newcommand{\lhat}{\mathcal F{L^{1}}}

\newcommand{\lloc}{\mathcal FL^{1}_{loc}}

\newcommand{\g}{\gamma}

\newcommand{\T}{\mathcal{T}}

\newcommand{\tB}{\widetilde{B}}
\newcommand{\tA}{\widetilde{A}}
\newcommand{\tE}{\widetilde{E}}

\DeclareMathOperator{\dist}{dist }

\DeclareMathOperator{\supp}{supp}

\begin{document}

%
%
%
%
%
%
%
%
%


\title [Closed Operator Functional Calculus]
{
Closed operator functional calculus \\   in Banach modules and applications} 

\author{Anatoly G. Baskakov}
\address{Department of Applied Mathematics and Mechanics \\ Voronezh State University \\ Voronezh 394693 \\ Russia}
\email{anatbaskakov@yandex.ru}



\author{Ilya A. Krishtal}
\address{Department of Mathematical Sciences\\ Northern Illinois University\\ DeKalb, IL 60115 \\ USA}
\email{ikrishtal@niu.edu}
\author{Natalia B. Uskova}
\address{Department of Higher Mathematics and Mathematical Physical Modeling \\ Voronezh State Technical University \\ Voronezh 394026 \\ Russia}
\email{nat-uskova@mail.ru}


\date{\today }




\begin{abstract}
 We describe a closed operator functional calculus in Banach modules over the group algebra $L^1(\RR)$ and  illustrate its usefulness with a few applications. In particular, we deduce a spectral mapping theorem for operators in the functional calculus, which generalizes some of the known results. We also obtain an estimate for the spectrum of a perturbed differential operator in a certain class.
\end{abstract}
\maketitle

\section{Introduction}\label{intro}

The goal of this paper is to describe a closed operator functional calculus in Banach modules over the group algebra $L^1(\RR)$ and to illustrate its usefulness with a few applications. The functional calculus was introduced in \cite{BK14} in order to obtain several non-commutative extensions of Wiener's $1/f$ lemma \cite{W32}.
The first application discussed in this paper (see Theorems \ref{smt1} and \ref{sad}) gives several versions of the spectral mapping theorem \cite{B79, B04, EN00, LMF73}.

The second application of the functional calculus is an estimate of the spectrum $\s(\mathscr L)$ of a differential operator   $\mathscr L = A - V: D(A)\subseteq L^2(\RR) \to L^2(\RR)$, $A = -i\frac d{dt}$. The domain $D(A)$ is chosen to be the Sobolev space $W^{1,2}(\RR)$ of absolutely continuous functions with the (almost everywhere) derivative in $L^2(\RR)$, and the operator $V: D(A)\subseteq L^2(\RR) \to L^2(\RR)$ is chosen to be of the form
\begeq\label{opv}
(Vx)(t) = v(t)x(-t), \ t\in\RR, v\in L^2(\RR).
\eq

In \cite{BKR17, BKU18}, we performed the spectral analysis for an analogous operator on $L^2([0,\omega])$. In fact, in that case, $\s(\mathscr L)$  is discrete and differs from $\sigma(A)$  by an $\ell^2$ sequence. Here, $\sigma(A) = \RR$, and we end up estimating a region in $\CC$ that contains $\s(\mathscr L)$. We cite \cite{BKR17, BKU18, BKh11z, KS17, SDI18} and references therein for the motivation of studying differential operators with an involution, such as the reflexion operator $V$.

The resulting estimate is contained in the following theorem.

\bt\label{maint}
Consider the operator
$\mathscr L = -i\frac d{dt} - V: W^{1,2}(\RR)\subseteq L^2(\RR) \to L^2(\RR)$ with $V$ of the form \eqref{opv}.
Then there exists a continuous real-valued function $f\in L^2(\RR)$ such that for any $\l\in\s(\mathscr L)$ one has
$
|\Im m\l| \le f(\Re e\l).
$
\et

Thus, in Theorem \ref{maint}, the spectrum $\s(\mathscr L)$ lies between the graphs of the functions $f$ and $-f$.

The remainder of the paper is organized as follows. In Section \ref{calc},  we introduce the necessary notions and notation and describe the functional calculus. In Section \ref{spacmap}, we formulate and prove  a few novel versions of the spectral mapping theorem. In Section \ref{simtrans}, we prove Theorem \ref{maint}. Finally, Section \ref{app} contains proofs of a few auxiliary results.

\section{$\mathcal F{L^1_{loc}}$ functional calculus.}\label{calc}

In our exposition of the closed operator functional calculus for generators of Banach $L^1(\RR)$-modules, we follow \cite{BK14}. 
Let us introduce some  notation.

We denote by $\X$ a complex Banach space and by $B(\X)$ the Banach algebra of all bounded linear operators in $\X$. We also assume that $\X$ is endowed with a non-degenerate Banach module structure over the group algebra $L^1(\RR)$. The multiplication in $L^1(\RR)$ is the convolution 
\[
(f*g)(t) = \int_\RR f(s)g(t-s)ds, \ f,g\in L^1(\RR),\ t\in\RR.
\]
\bd\label{dmod}
A complex Banach space $\X$  is a Banach module over $L^1(\RR)$ if there is a bilinear map
$(f,x)\mapsto fx: L^1(\RR)\times\mathcal X\to\mathcal X$ which has the following properties:
\begin{enumerate}
\item $(f*g)x = f(gx)$,  $f,g\in L^1(\RR)$, $x\in\mathcal X$;
\item $\norm{fx} \le \|f\|_1\|x\|$,  $f\in L^1(\RR)$, $x\in\mathcal X$.
\end{enumerate}
\ed

As usual (see \cite{BD19, BK05} and references therein), by 
non-degeneracy of the module we mean that $fx = 0$ for all $f\in L^1(\RR)$ implies that $x = 0$. We only consider Banach module structures that are associated with an  isometric representation $\T: \RR\to B(\X)$, that is we have
\begeq\label{assrep}
\T(t)(fx) = f_t x = f(\T(t)x),\ t\in\RR, f\in L^1(\RR), x\in\X,
\eq
where $f_t(s) = f(t+s)$, $t,s\in\RR$. With a slight abuse of notation \cite{BK05}, given $f\in L^1(\RR)$, we shall denote by $\T(f)$ the operator in $B(\X)$ defined by $\T(f)x = fx$, $x\in\X$. Observe that we have $\|\T(f)\|\le \|f\|_1$, $f\in L^1(\RR)$, by Property 2 in Definition \ref{dmod}. For the Banach module $\X$, we will also use the notation $(\X,\T)$ if we want to emphasize that the module structure is associated with the representation $\T$.

We use the Fourier transform of the form
\[
(\F(f))(\xi) = \widehat f(\xi) = \int_\RR f(t)e^{-it\xi}dt, \ f\in L^1(\RR),
\]
so that 
$\|\widehat f\|_2 = \sqrt{2\pi}\|f\|_2$, $f\in L^2(\RR)$. We shall denote by $\mathcal FL^1 = \mathcal FL^1(\RR)$ the Fourier algebra $\mathcal F(L^1(\RR))$. The inverse Fourier transform of a function $h\in\lhat(\RR)$ will be denoted by $\check h$ or $\F^{-1}(h)$.

\bd
Let $(\X,\T)$ be a non-degenerate Banach $L^1(\RR)$-module, and $N$ be a subset of $\X$. The \emph{Beurling spectrum} $\Lambda(N) = \Lambda(N,\T)$ is defined by
\[
\Lambda(N,\T) = \{\l\in\RR: fx = 0 \mbox{ for all } x\in N \mbox{ implies } \widehat f(\l) = 0, f\in L^1\}.
\] 
\ed


To simplify the notation we shall write $\Lambda(x)$ instead of $\Lambda(\{x\})$, $x\in\X$. 
We refer to \cite[Lemma 3.3]{BK05} for the basic properties of the Beurling spectrum.
We  also define
\[
\X_{comp} = \{x\in\X: \Lambda(x) \mbox{ is compact}\},\ \X_\Phi = \{\T(f)x: f\in L^1(\RR), x\in \X\}
\]
 and
\[
\X_c = \{x\in \X: \mbox{ the function } t\mapsto \T(t)x: \RR\to\X \mbox{ is continuous}\},
\] 


For any $z\in\CC\setminus\RR$,  consider the function $f_z\in L^1(\RR)$ whose Fourier transform is the function $\phi_z: \RR\to\CC$ defined by $\phi_z(\lambda) = (\lambda - z)^{-1}$, $\lambda\in\RR$. Hilbert's resolvent identity holds for the operator-valued function
$R: \CC \setminus \RR\to B(\X)$ given by $R(z) = \T(f_z)$, $z\in\CC \setminus \RR$. Since the $L^1(\RR)$-module $\X$ is non-degenerate, we have $\bigcap_{z\in\CC\setminus\RR}\ker R(z) = \{0\}$.
Therefore \cite{B13}, $R$ is the resolvent of some linear operator $\A: D(\A)\subseteq \X\to \X$. This operator $\A$ is called the generator of the $L^1(\RR)$-module $\X$. 
We remark that if $\T: \RR\to B(\X)$ is a strongly continuous group representation, then $i\A$ is its generator.


It is not hard to show that the operators $\T(f)$, $f\in L^1(\RR)$, provide a functional calculus for the generator $\A$. Via the isomorphism of $L^1(\RR)$ and $\mathcal FL^1(\RR)$, we also get the functional calculus $\chT(\widehat f) = \T(f)$, $\widehat f = \mathcal F(f)\in\mathcal FL^1$. It is useful to extend this functional calculus to the space $\lloc(\RR)=\{h$: $\RR \to \CC$ such that  
$h\widehat\varphi\in\lhat(\RR)$  for any $\varphi\in L^1(\RR)$
 with 
$\supp{\widehat\varphi}$ compact$\}$. Observe that $\lhat(\RR)\subset\lloc(\RR)$. Moreover,
$\lloc$ is also an algebra under pointwise multiplication.

For $h\in\lloc(\RR)$ we define a (closed) operator $\chT(h) = h\diamond : D(h\diamond)=D(\chT(h))\subseteq\X\to\X$ in the following way.
First, let $x\in\X_{comp}$ and 
\begeq\label{cl}
\chT(h)x = h\diamond x : = (h\widehat\varphi)^{\vee} x = \T((h\widehat\varphi)^{\vee})x,
\eq
where $\varphi\in L^1(\RR)$ is such that $\supp \widehat\varphi$ is compact and  $\widehat\varphi\equiv 1$ in a neighborhood of $\Lambda(x)$. The vector $\chT(h)x$ is well defined in this way because it is independent of the choice of $\varphi$. 

Next, we extend the definition of $\chT(h)$ by taking the closure of the just defined operator on $\X_{comp}$. In other words,
if $x_n\in \X_{comp}$, $n\in\NN$, $x = \lim\limits_{n\to \infty} x_n$, and
$y = \lim\limits_{n\to \infty} h\diamond x_n$ exists, we let $\chT(h)x = h\diamond x = y$. 
Lemma 2.7 in \cite{BK14} shows that $\chT(h)$ is then a well-defined closed linear operator and we do, indeed, have $\chT(\widehat f) = \T(f)$, $f\in L^1(\RR)$. Moreover, applying \cite[Proposition 2.8]{BK14}, we get for a given $x\in D(\chT(h))$ that
$\T(t)(h\diamond x) = h\diamond (\T(t)x))$ and
\begeq\label{func}
\T(f)(h\diamond x) = h\diamond (\T(f)x))=(\widehat f h)\diamond x,
\eq
$t\in\RR$, $f\in L^1(\RR)$, $h\in\lloc(\RR)$. We note that $D(\chT(h))\subseteq \X_c$ for all $h\in \lloc$.
We also note the following useful property that is implied by the definition of the operators $\chT(h)$ and \cite[Lemma 3.3]{BK05}:
\begeq\label{hex}
\L(\chT(h)x, \T)\subseteq \supp h\cap\L(x,\T),\ h\in \lloc, x\in\X.
\eq

It is  not hard to see 
that the generator $\A$ of the module $(\X,\T)$ satisfies 
\begeq\label{ida}
\A = \chT(\id),
\eq 
where $\id\in\lloc(\RR)$ is the identity function $\id(\xi) = \xi$, $\xi\in\RR$. Thus, we have an $\lloc$ functional calculus for the generator $\A$, which we will use to prove a few spectral mapping results and construct a similarity transform to obtain an estimate for the spectrum of the perturbed differential operators.

We will use the following  sufficient condition for functions in $L^2$ to belong to $\lhat$.
For completeness, we provide its proof in Section \ref{app}.

\bl\label{funest}
Assume $f\in L^2(\RR)$ and $\widehat f \in W^{1,2}(\RR)$. Then $f\in L^1(\RR)$ and
\begeq\label{estfun}
\|f\|_1^2 \le 2\|\widehat f\|_2\|\widehat f'\|_2.
\eq
\el

We  illustrate the above lemma with the following two examples.

\bex\label{trapez}
For $a > 0$,
consider the ``trapezoid function'' $\tau_a$ defined by
\[
\tau_a(\xi) = 
\begin{cases}
1, & |\xi|\le a, \\
\frac1a(2a-|\xi|), & a<|\xi|< 2a,\\
0, & |\xi|\ge 2a.
\end{cases}
\]
Direct computations show that $\|\tau_a\|_2 = 2\sqrt{\frac23a}$, $\|\tau'_a\|_2 = \sqrt \frac2a$, and $\tau_a = \widehat \varphi_a$, where
\[
\varphi_a(t) = \frac{2\sin \frac{3 at}2  \sin\frac {at}2}{\pi at^2}, \ t\in\RR.
\]
From Lemma \ref{funest} we conclude that $\|\varphi_a\|_1 \leq 2^{\frac32}\cdot 3^{-\frac14}$. We remark that 
\cite[Lemma 1.10.1]{Le53} yields a better estimate: $\|\varphi_a\|_1 \leq \frac4\pi+\frac2\pi\ln 3$. An even better estimate, $\|\varphi_a\|_1 \leq \sqrt 3$, follows from
\cite[Proposition 5.1.5]{RS00}. 
\eex

\bex\label{ex2}
For $a > 0$, let 
\[
\omega_a(\xi) = \frac1\xi(1-\tau_a(\xi))= 
\begin{cases}
0, & |\xi|\le a, \\
-\frac1a-\frac1{\xi}, & -2a<\xi\le -a,\\
\frac1a-\frac1{\xi}, & a<\xi\le 2a,\\
\frac1{\xi}, & |\xi|> 2a.
\end{cases} 
\]
Then $\|\omega_a\|_2 = \sqrt{\frac{4-4\ln 2}a}\le 1.11/\sqrt a$ and $\|\omega'_a\|_2 = \sqrt{\frac{2}{3a^3}}\le 0.82/(a\sqrt a)$. It follows that the functions $\psi_a$ defined by $\widehat\psi_a = \omega_a$ satisfy
$\|\psi_a\|_1^2\le \frac 4{a^2}\sqrt{\frac23(1-\ln 2)}$ so that 
$ \|\psi_a\|_1\le 1.35/a.$

We will also need an estimate for $\|\psi_a\|_\infty$. Observe that for $t>0$
\[
\bs
\frac1{2\pi}&\left|\int_{-\infty}^{\infty} \omega_a(\xi)e^{it\xi}d\xi\right|  = 
\frac1\pi \left|\int_{0}^{\infty} \omega_a(\xi)\sin(t\xi)d\xi\right|  \le \\
\frac1\pi &\left(1+\left|\int_{a}^{\infty} \frac{\sin(t\xi)}{\xi}d\xi\right|\right) =
\frac1\pi \left(1+\left|\int_{at}^{\infty} \frac{\sin\xi}{\xi}d\xi\right|\right)\le \frac1\pi+1.
\end{split}
\]
Since $\omega_a$ is an odd function, it follows that
\begeq
 \|\psi_a\|_\infty\le \frac1\pi+1.
\eq
\eex

Now we use the functions from the above two examples in our $\lloc$ functional calculus. In view of \eqref{func} and \eqref{ida}, we get the following crucial relationship:
\begeq\label{eqmet}
\bs
\A\T(\psi_a)x & 
= \chT(\id)\T(\psi_a)x  
= \chT(\id \cdot\omega_a)x = \chT(\mathds 1-\tau_a)x = x - \T(\varphi_a)x,
\end{split}
\eq
which holds for every $x\in D(\A)$; by $\mathds 1\in \lloc$ we denoted the function $\mathds 1(\xi) = 1$, $\xi\in \RR$. In fact, since $D(\chT(\mathds 1)) = \X_c$ and $\A$ is a closed operator, we get that \eqref{eqmet} holds for all $x\in\X_c$. Moreover, non-degeneracy of the module $(\X, \T)$ implies that the  operator
$\A\T(\psi_a)$ extends uniquely to a bounded operator $I-\T(\varphi_a)\in B(\X)$. To simplify the notation, given $\l\in\CC$, we may write $\l - f$ instead of $\l\eye-f$ for functions and $\l-A$ instead of $\l I-A$ for operators. 

We note that the family $(\varphi_\a)$ from Example \ref{trapez} has another useful property: it forms a bounded approximate identity.

\bd\label{bai}
A family of functions $(\phi_a)_{a > 0}$ is called a \emph{bounded approximate identity} or b.a.i.~if $\|\varphi_a\|_1 \le M$ for all $a > 0$ and $\lim_{a\to \infty} \|\varphi_a*f - f\|_1 =0$ for all $f\in L^1(\RR)$. 
\ed

Following \cite{BK14}, we call a \bai $(\phi_\a)$ a \cfbai if $\supp \widehat\phi_a$ is compact for each $a>0$. Clearly, the family $(\varphi_\a)$ from Example \ref{trapez} is a \cfbai Another useful \cfbai is given by the family
\begeq\label{tri}
\g_a(t) =\frac a{2\pi}\left(
\frac{\sin(at/2)}{at/2}\right)^2,\ t\in\RR, a>0.
\eq
Observe that $\widehat\g_a(\xi) = (1-|\xi|/a)\eye_{[-a,a]}(\xi) =: \triangle_a(\xi)$ is the so called triangle function and, since $\g_a$ is non-negative, $\|\g_a\|_1 = \triangle_a(0) = 1$.

We conclude this section by recalling the following  result that contains the celebrated Cohen-Hewitt factorization theorem \cite{C59, H64}.

\bp[\cite{BK05}, Lemma 4.3]\label{cht}
For any \bai $(\phi_a)$, we have
\[
\X_c = \X_\Phi = \overline{\X_{comp}} = \{x\in \X: x =\lim_{a\to\infty} \phi_a x = x\}.
\]
\ep

\section{Spectral mapping theorems.}\label{spacmap}

We begin this section by recalling  a spectral mapping theorem from \cite{B04} which we endeavor to extend (see also \cite{EN00, LMF73}). 

\bt[\cite{B04}, Corollary 1.5.3]\label{smt0}
Let  $\X$ be a non-degenerate Banach $L^1(\RR)$-module with the structure associated with a 
representation $\T$.
For $f\in L^1(\RR)$,   $\s(\T(f)) = \overline{\widehat f(\L(\X,\T))}$. 
\et

To prove our extensions we will need the following two lemmas that, in particular, give a special case of the above result. 

\bl\label{sml1}
Assume that $K = \L(\X,\T)$ is compact. Then,
for $h\in \lloc$, we have $\s(\chT(h))\subseteq  h(K) $.
\el

\bpf
Observe that $h(K)$ is automatically compact as a continuous image of a compact set.
Assume $\l\notin h(K)$. Then the function $u$ given by $u(t) = \frac1{\l-t}$ is analytic in a neighborhood of the compact set $h(K)$. The Wiener-L\'evy theorem \cite[Theorem 1.3.1]{RS00} ensures existence of a function $f\in L^1(\RR)$ such that $\widehat f(\xi) = \frac1{\l-h(\xi)}$ for $\xi$ in a neighborhood of $K$. For $x\in\X =\X_{comp}$, we use   \eqref{func} and \cite[Lemma 3.3]{BK05} to obtain
\[
(\l I-\T(f))\chT(h)x = \chT(h)(\l I-\T(f))x = \chT(\l\widehat f - \widehat f h)x = \chT(\eye)x = x,
\]
 so that $\l\in\rho(\chT(h))$. Thus, $\s(\chT(h))\subseteq  h(K)$.
\epf

\bl\label{sml2}
For $h\in \lloc(\RR)$, we have $\s(\chT(h))\supseteq \overline{h(\L(\X,\T))}$.
\el

\bpf
Assume $\l =  h(\xi) \in\rho(\chT(h))$ for some $\xi\in\RR$ and let $f\in L^1(\RR)$ be such that $\supp \widehat f$ is compact, $\xi\in\supp \widehat f$, and ${\widehat f(\supp  h)}\subset \rho(\chT(h))$. Let also $\Y = \T(f)\X$ and $B = \chT(h)\vert_\Y$ be the restriction of the operator $\chT(h)$ to the submodule $\Y$. Clearly, $B\in B(\Y)$. We claim that $\rho(B) = \CC$, which would imply $\Y = \{0\}$ yielding $\supp \widehat f\cap \L(\X,\T) = \emptyset$. To prove the claim, we first observe that $\rho(\chT(h))\subseteq\rho(B)$. Indeed, since $\chT(h)$ commutes with $\T(f)$ by \eqref{func}, we have that for $\l\in\rho(\chT(h))$ the resolvent operator $(\l I - \chT(h))^{-1}$ also commutes with $\T(f)$ ensuring $(\l I - B)^{-1} = (\l I - \chT(h))^{-1}\vert_\Y$. Using \eqref{hex}, we get 
\[
{\widehat f(\L(\Y))} \subseteq {\widehat f(\supp  h)}\subset \rho(\chT(h)) \subseteq \rho(B). 
\]
Secondly, Lemma \ref{sml1} implies ${\widehat f(\L(\Y))}^c \subseteq \rho(B)$, and the claim is established. It follows that $\supp \widehat f\cap \L(\X,\T) = \emptyset$, and hence $\xi\notin \L(\X,\T)$. Thus, $\s(\chT(h))\supseteq {h(\L(\X,\T))}$ and, since the spectrum is closed, $\s(\chT(h))\supseteq \overline{h(\L(\X,\T))}$.
\epf

%

Thus, to extend Theorem \ref{smt1} to the $\lloc$ setting we only need an
 analog of Lemma \ref{sml1} for the case when $\L(\X)$ is not necessarily compact. This, however, may not always hold at this level of generality as we can no longer use the Wiener-L\'evy theorem. We offer several ways to circumvent the problem.
 
 First, we present a result that is immediate from the proof of Lemma \ref{sml1}.
 
 \bp\label{smp1}
 Let $h\in \lloc$ and $\l\in\CC$. Assume that  there exists a function $g_\l\in\lloc$ such that $g_\l = (\l-h)^{-1}$ in a neighborhood of $\L(\X, \T)$ and $\chT(g_\l)$ belongs to $B(\X)$. Then $\l\in\rho(\chT(h))$.
 \ep
 

This motivates the following definition.

\bd
Let $\X =(\X,\T)$ be a non-degenerate Banach $L^1(\RR)$-module. A function $h\in\lloc(\RR)$ is called
\emph{$\X$-regular} if for any $\lambda \notin \overline{h(\L(\X,\T)}$ there exists a function $g_\l\in\lloc(\RR)$ such that 
\begeq\label{gl}
g_\l(\xi)(\l-h(\xi))=1, 
\eq   
for every $\xi$ in a neighborhood of $\L(\X,\T)$, and $\chT(g_\l)$ belongs to $B(\X)$.
\ed

Clearly, it would be sufficient for the functions $g_\l$ in the above definition to belong to $\lhat$. Hence, by the Wiener-L\'evy theorem, if $\L(\X)$ is compact, every $h\in\lloc(\RR)$ is 
{$\X$-regular}. 
%
%

The next result is now immediate.

\bt\label{smt1}
Assume that $h\in\lloc(\RR)$ is $\X$-regular. Then $\s(\chT(h))= \overline{h(\L(\X,\T))}$. Moreover, given $\lambda \notin \overline{h(\L(\X,\T))}$, we have $(\l-\chT(h))^{-1} = \chT(g_\l) = \T(\check g_\l)$, where $g_\l\in \lloc(\RR)$ is defined by \eqref{gl}.
\et

The following result shows that Theorem \ref{smt1} does indeed generalize Theorem \ref{smt0}.

\bp
Any function $h\in \lhat(\RR)$ is $\X$-regular for any $\X$. 
\ep

\bpf
Observe that by the Riemann-Lebesgue lemma  if $0\notin\overline{h(\L(\X,\T))}$ then $\L(\X,\T)$ is compact, and the result follows.

Assume now that $0\neq \l\notin \overline{h(\L(\X,\T))}$. Then, without loss of generality we may assume that $\l\notin \overline{h(\RR)}$. Indeed, if that was not the case, we would have $0\neq \l \in h(\RR)$ and $h^{-1}(\{\l\})$ would be a compact set disjoint from $\L(\X,\T)$. We could then find $\phi \in\lhat$ with compact support that is disjoint from $\L(\X,\T)$ and such that $\l\notin \overline{(h+\phi)(\RR)}$. We would then apply the following argument to $h+\phi$ instead of $h$.

A modification of the Wiener-L\'evy theorem (see \cite[Theorem 1.3.4]{RS00} or \cite{GRS64}) or a special case of the Bochner-Phillips theorem (see \cite[Theorem 10.3]{BK14} or \cite{BP42}) show that 
\[
g_\l = (\l-h)^{-1} = \l^{-1} +\widetilde h
\]
for some $\widetilde h \in\lhat$. Then $\chT(g_\l) = \l^{-1}I+\chT(\widetilde h) \in B(\X)$, and the result follows.
\epf

Another  sufficient condition for $\X$-regularity follows from Lemma \ref{funest}.

\bp\label{funp}
Assume that $h\in\lloc(\RR)$ is such that for every $\lambda \notin \overline{h(\L(\X,\T))}$ there exists a function $g_\l\in W^{1,2}(\RR)$ that satisfies \eqref{gl} in a neighborhood of $\overline{h(\L(\X,\T))}$. Then $h$ is $\X$-regular. In particular, every polynomial is $\X$-regular  for any $\X$. 
\ep

\bpf
The first assertion follows immediately from Lemma \ref{funest}. To prove the second one, we 
note that for every polynomial $p$ and $\lambda \notin \overline{p(\RR)} = p(\RR)$ the function $(\l-p(\cdot))^{-1}$ belongs to $W^{1,2}(\RR)$. This shows that $p$ is $\X$-regular  for any $\X$ such that $\L(\X,\T) = \RR$.
In case $\L(\X,\T) \neq \RR$, given 
$\l \in {p(\RR)}\setminus \overline{p(\L(\X,\T))}$, we observe that $p^{-1}(\{\l\})$ is a finite set that does not intersect $\L(\X,\T)$. Hence, there is an infinitely many times differentiable 
function 
$\phi _\l\in\lhat$ with compact support that is disjoint from $\L(\X,\T)$ and such that $\l\notin \overline{(p+\phi_\l)(\RR)}$.
Then $g_\l = (\l-p-\phi_\l)^{-1}\in W^{1,2}(\RR)$, and  the result follows.
\epf

The following well-known result now follows immediately from \eqref{ida} and the fact that the Beurling spectrum is a closed set.

\bc
The generator $\A$ of a non-degenerate Banach $L^1(\RR)$-module $\X$ satisfies
$
\s(\A) = {\L(\X)}.
$
\ec

In the context of Proposition \ref{funp}, Lemma \ref{funest}  also allows us to estimate the resolvent of the operators $\chT(h)$.

\bc
Assume that $h\in\lloc(\RR)$, $\lambda \notin \overline{h(\RR)}$, and the function $g_\l$ defined by \eqref{gl} belongs to $W^{1,2}(\RR)$. Then
\[
\left\|(\l-\chT(h))^{-1}\right\| \le\sqrt{2\|g_\l\|_2\|g'_\l\|_2}. 
\]
\ec

The following definition allows us to provide yet another example of $\X$-regularity.

\bd We say that $h\in\lloc(\RR)$ is an \emph{almost periodic function with a summable Fourier series} if
\begeq\label{aph}
h(\xi) = \sum_{n\in\ZZ} c_n e^{i\xi t_n}, \ \sum_{n\in\ZZ} |c_n| < \infty,\ \xi, t_n\in\RR, n\in\ZZ.
\eq
The set of all such functions is denoted by $\AP_1$ or $\AP_1(\RR)$. 
\ed

We note that $\AP_1$ is a Banach space with the norm 
\[
\|h\|_{\AP_1} = \sum_{n\in\ZZ} |c_n|,
\]
where $h\in \AP_1$ is given by \eqref{aph}. We also mention 
 \cite[Proposition 2.11]{BK14}, which states that for such $h$ we have
\begeq\label{mas}
\chT(h) = \sum_{n\in\ZZ} c_n\T(t_n) \in B(\X).
\eq

\bp\label{csp1}
Any function $h\in \AP_1$ is $\X$-regular for any $\X$. Moreover, if $\l\notin\overline{h(\RR)}$, then 
\[
\norm{\left(\l-\chT(h)\right)^{-1}} \le \norm{\frac1{\l-h}}_{\AP_1}.
\]
\ep

\bpf
Let $\RR_d$ be the group of real numbers with the discrete topology and $\RR_c$ be its Pontryagin dual -- the Bohr compactification of $\RR$.
It is well-known that a function in $\AP_1(\RR)$ has a unique continuous extension to  $\RR_c$ and can be identified with an element of $\lhat(\RR_d)$ -- the Fourier algebra of the group $\RR_d$.   The closure of $\L(\X, \T)$ in $\RR_c$ is then a compact subset of $\RR_c$ and, given $\l\notin \overline{h(\L(\X,\T))}$, 
the Wiener-L\'evy theorem for locally compact Abelian groups (\cite[Theorem 6.1.1]{RS00}) establishes existence of $g_\l \in\AP_1(\RR)$ that satisfies \eqref{gl} in a neighborhood of $\L(\X,\T)$ in $\RR$. Hence, $\X$-regularity follows from  \eqref{mas}, i.e.~\cite[Proposition 2.11]{BK14}.

For $\l\notin\overline{h(\RR)}$, it suffices to apply 
 the almost periodic version of Wiener's $1/f$ lemma \cite{BK10, L53} that shows that $g_\l = \frac1{\l-h}\in \AP_1$. The desired estimate then follows from 
 \eqref{mas}.
\epf


\bex
Let $h(\xi) = e^{i\xi t_0}$ for some $t_0 \in \RR$ and $\l\in\CC\setminus\TT$, where $\TT = \{z\in\CC: |z|=1\}$. From \eqref{mas}, we get $\chT(h) = \T(t_0)$. Using the estimate from Proposition \ref{csp1}, we get by direct computation that 
\[
\norm{\left(\l-\T(t_0)\right)^{-1}} \le\left\|\frac{1}{\l-h}\right\|_{\AP_1} = \frac{1}{|1-|\l||} = \left(\dist(\l,\TT)\right)^{-1}.
\]
\eex

In general, it may be hard to check if a given function is $\X$-regular.  We cite \cite{BD19} and references therein for related results. We note that the notion of regularity at infinity discussed in \cite{BD19} is more restrictive than $\X$-regularity.

The following theorem provides a special case when the assumption of $\X$-regularity is not needed.

\bt\label{smt2}
Assume that $\X = H$ is a Hilbert space and the representation $\T$ is unitary. Then for any $h\in\lloc(\RR)$ we have $\s(\chT(h))= \overline{h(\L(\X,\T))}$. Moreover,  given $\lambda \notin \overline{h(\L(\X,\T))}$, we have
\begeq\label{be1}
\left \|(\l-\chT(h))^{-1}\right\| = \left(\dist\left(\l, \overline{h(\L(\X,\T)}\right)\right)^{-1}.
\eq
\et

\bpf
In view of Lemma \ref{sml2}, we only need to prove  $\s(\chT(h))\subseteq \overline{h(\L(\X,\T))}$, $h\in \lloc(\RR)$. Pick $\l\notin \overline{h(\L(\X,\T))}$. We will show that $\l\in\rho(\chT(h))$.

Let $(\phi_a)$ be a \cfbai  and 
$\X_a = \T(\phi_a)\X$, $a > 0$, be the corresponding submodules of $\X$. From 
\eqref{hex}, we have $\L(\X_a, \T)\subseteq \supp\widehat \phi_a\cap \Lambda(\X,\T)$. Hence, $\X_a\subseteq\X_{comp}\subseteq D(\chT(h))$ and \eqref{func} implies that $\X_a$ is invariant for $\chT(h)$. Therefore, the restrictions of $\chT(h)$ to $\X_a$, $a> 0$, are well defined. We will denote these restrictions by $B_a$. 

Since $h$ is $\X_a$-regular,
Theorem \ref{smt1} applies for $B_a$ yielding 
$\s(B_a)= h(\L(\X_a, T))\subseteq \overline{h(\L(\X,\T))}$. It follows that
$\lambda\in \rho(B_a)$.
Moreover, since the representation $\T$ is unitary, the operators $B_a$ are normal. Therefore,
the norms of their resolvents satisfy
\begeq\label{resest}
\|R(\l; B_a)\| = (\dist(\l, \s(B_a))^{-1} \le \left(\dist\left(\l, \overline{h(\L(\X,\T)}\right)\right)^{-1}.
\eq
Now, since the representation $\T$ is strongly continuous, 
Proposition \ref{cht} implies that 
\[
x = \lim_{a\to\infty} \T(\phi_a)x
\]
for an arbitrary $x\in \X$. From \eqref{resest} and the Banach-Steinhaus theorem, we get that $C = C_\l$ given by
\begeq\label{strl}
Cx = \lim_{a\to\infty}  R(\l; B_a)\T(\phi_a)x
\eq
is a well-defined bounded linear operator. By direct computation, it follows that
\[\T(f)C(\l-\chT(h))x = \T(f)(\l-\chT(h))Cx, 
\]
for any $x\in D(\chT(h))$ and $f\in L^1(\RR)$ with $\supp\widehat f$ compact. Since the module $\X$ is non-degenerate, we get $C = (\l-\chT(h))^{-1}$. 
Finally, the equality in \eqref{be1} follows since $C=(\l-\chT(h))^{-1}$ is a normal operator.
\epf

\brem
Often \cite{ABK08, B97Sib, G04} a representation $\T$ and operators of the form $\chT(h)$ act not just in a single Banach module but in a whole chain $(\X_p)$ of such modules. For example, matrices with sufficient off-diagonal decay  define bounded operators on all $\ell^p$, $p\in[1,\infty)$. In this case, it is not unusual for $\s(\chT(h))$ to be independent of $p$. If also one of the modules $\X_p$ happened to be a Hilbert space, Theorem \ref{smt2} would then yield a spectral mapping theorem for all Banach modules $\X_p$ in the chain.
\erem

The proof of Theorem \ref{smt2} leads us to define the following notion of regularity for functions in $\lloc(\RR)$.

\bd\label{speca}
Let $(\varphi_a)$ be the \cfbai from Example \ref{trapez}.
A function $h\in\lloc(\RR)$ is called \emph{spectrally admissible} if 
for each $\lambda \notin \overline{h(\L(\X,\T))}$ there exist  functions $g_\l^a\in\lhat(\RR)$ such that $g_\l^a= \frac{\widehat\varphi_a}{\l-h}$ in a neighborhood of  $\overline{h(\L(\X,\T))}$ and
\begeq\label{mhl}
M_h(\l): = \sup_{a > 0} \left\|\F^{-1}\!\left(g_\l^a\right)\right\|_1 < \infty.
\eq 
\ed

\brem\label{anybai}
In the above definition, instead of the functions from Example \ref{trapez} we may use
 a \cfbai $(\varphi_{a,n})$, $n>1$, given by
\[
\widehat \varphi_{a,n} (\xi)= \tau_{a,n}(\xi) = 
\begin{cases}
1, & |\xi|\le a, \\
\frac1{(n-1)a}(na-|\xi|), & a<|\xi|< na,\\
0, & |\xi|\ge na.
\end{cases}
\]
From
\cite[Proposition 5.1.5]{RS00} we get  $\|\varphi_{a,n}\|_1 \leq \sqrt {\frac {n+1}{n-1}}$,
which may give a  smaller $M_h(\l)$ in \eqref{mhl}.
We also note that  Lemma \ref{funest} may often be used to prove spectral admissibility.
\erem

\bt\label{sad}
Let $(\X,\T)$ be a non-degenerate Banach $L^1(\RR)$-module such that the representation $\T$ is strongly continuous.
Assume that   a function $h\in\lloc(\RR)$ is {spectrally admissible}. Then
$\s(\chT(h))= \overline{h(\L(\X,\T))}$. Moreover, given $\lambda \notin \overline{h(\L(\X,\T))}$, we have
\[
\left \|(\l-\chT(h))^{-1}\right\| \le M_h(\l).
\]
\et
\bpf
As in the proof of Theorem \ref{smt2}, given a \cfbai $(\phi_a)$, we let $\X_a = \T(\phi_a)\X$ and $B_a = \chT(h)|_{\X_a}$. 
Observe that for a sufficiently large $b > 0$, the function $g_\l^b$ from Definition \ref{speca} satisfies $g_\l^b (\xi)= \frac1{\l-h(\xi)}$ for every $\xi$ in a neighborhood of $\overline{h(\L(\X_a,\T))}$.
From
Theorem \ref{smt1}, we deduce that $$\|R(\l, B_a)\| \le \left\|\F^{-1}\!\left(g_\l^{b}\right)\right\|_1 \le M_h(\l) < \infty,\ \lambda \notin \overline{h(\L(\X,\T))}.$$
The remainder of the proof of  Theorem \ref{smt2} now goes through in this setting. An application of the Banach-Steinhaus theorem  shows that  \eqref{strl} defines an operator $C_\l\in B(H)$ satisfying 
$\|C_\l\| \le M_h(\l)\sup_a\|\phi_a\|$, and, choosing $\phi_a = \g_a$ defined by \eqref{tri} gives $\|C_\l\| \le M_h(\l)$. It is then verified by direct computation that $C_\l = (\l-\chT(h))^{-1}$.
\epf

\section{Spectral estimates for the operator $\mathscr L$.}\label{simtrans}

In this section, we prove Theorem \ref{maint}.
The approach we pursue is based on the following result which holds for Hilbert-Schmidt perturbations of  general self-adjoint operators on an abstract complex Hilbert space $H$. 
The ideal of all Hilbert-Schmidt operators in $H$ will be denoted by $\mathfrak{S}_2(H)$.
\bt\label{hspert}
Let $A: D(A)\subseteq H\to H$ be a self-adjoint operator and $B\in \mathfrak{S}_2(H)$. Then there exists a continuous real-valued function $f\in L^2(\RR)$ such that for any $\l\in\s(A+B)$ one has
$
|\Im m\l| \le f(\Re e\l).
$
\et

We believe that the above result has been known for a long time. Since we didn't find the reference, however, we provide its proof in Section \ref{app}. We cite \cite{GK69} for related results.

Clearly, the perturbation $V$ of the form \eqref{opv} may not be Hilbert-Schmidt; it is not even a bounded operator, in general. We will, however, construct a similarity transform which will allow us to use the above result.

\bd 
Two linear operators $A_i: D(A_i) \subset H \to H$, $i = 1,2$, are called \emph{similar} if there exists an invertible operator $U \in B(H)$ such that $UD(A_2) = D(A_1)$ and $A_1Ux = UA_2x$, 
$x \in D(A_2)$. We call the operator $U$ the \emph{similarity transform} of $A_1$ into $A_2$.
\ed

It is immediate that for similar operators $A_1$ and $A_2$ one has $\s(A_1) = \s(A_2)$.
Thus, to prove Theorem \ref{maint}, it suffices to
construct a similarity transform of $\mathscr L$ into $-i\frac d{dt}- B$ with  
$B\in \mathfrak{S}_2(H)$. 
In order to do it, we apply the $\lloc$ functional calculus in the space $\X  = \mathfrak L_A(H)$ of closed linear $A$-bounded operators that is defined as follows.

\bd\label{abdd}
 Let $A: D(A) \subset H \to H$ be a closed linear operator. 
 A linear operator $X: D(X) \subset H \to H$ is $A$-bounded if $D(X) \supseteq D(A)$ and $\|X\|_A = \inf\{c > 0: \|Xx\| \le c(\|x\| + \|Ax\|),\ x \in D(A)\} < \infty$.
\ed
The space $\mathfrak L_A(H)$ of all $A$-bounded linear operators with the  domain equal to $D(A)$ is a Banach space with respect to the norm $\|\cdot\|_A$. For densely defined operators $A$, restricting the domain of bounded operators to $D(A)$, allows us to view $B(H)$ as a subspace of $\mathfrak L_A(H)$.

Now we need to define the Banach module structure in $\X = \mathfrak L_A(H)$ with $A = -i\frac d{dt}$. We begin with a Banach module structure in $H = L^2(\RR)$.

The operator $A = -i\frac d{dt}: W^{1,2}(\RR)\subseteq L^2(\RR)\to L^2(\RR)$  is self-adjoint, and the operator $iA$ generates an isometric strongly continuous group of translations $T: \RR\to B(H)$, $T(t)x(s) = x(t+s)$, $x\in L^2(\RR)$. 
The non-degenerate $L^1$-module structure in $H$ is then given by convolution:
\[
T(f)x = \int_\RR f(t)T(-t)xdt, \ x\in H.
\]
Next, we let $\T:\RR\to B(\X) = B(\mathfrak L_A(H))$ be defined by $\T(t)X = T(t)XT(-t)$, $t\in\RR$. Since $T$ is an isometric representation, we get that $\T$ also has this property. We then have that
\begeq\label{qq}
(\T(f)X)x = \int_\RR f(t)(\T(-t)X)xdt = \int_\RR f(t)T(-t)XT(t)xdt,\ 
\eq
$X\in\X, x\in D(A), f\in L^1(\RR),$ defines a non-degenerate $L^1$-module structure in $\X$ that is associated with the representation $\T$. Moreover, the generator $\A$ of the module $(\X,\T)$ satisfies
\[
\A X = AX - XA, \ X\in D(\A),
\]
see e.g. \cite{EN00}. We now apply \eqref{eqmet} in this Banach module $(\X, \T)$ to get
\begeq\label{eqmet1}
A(\T(\psi_a)X) - (\T(\psi_a)X)A = X - \T(\varphi_a)X,\ X\in D(\A),
\eq
where the functions $\varphi_a$ and $\psi_a$, $a >0$, are defined in Examples \ref{trapez} and \ref{ex2}.
Moreover, the discussion following \eqref{eqmet} shows that for any $X\in \mathfrak L_A(H)$ and $x\in D(A)$ we have
\begeq\label{eqmet2}
A(\T(\psi_a)X)x - (\T(\psi_a)X)Ax = Xx - (\T(\varphi_a)X)x.
\eq

\bl\label{ml}
Consider  the functions $\varphi_a$ and $\psi_a$, $a >0$,  defined in Examples \ref{trapez} and \ref{ex2}.
An operator $V$ of the form \eqref{opv} has the following properties.
\ben
\item $\T(\varphi_a)V \in \mathfrak S_2(H)$ and $\|\T(\varphi_a)V\|_2=\sqrt{\frac{2a}{3\pi}}\|v\|_2$.
\item $\T(\psi_a)V \in \mathfrak S_2(H)$ and $\|\T(\psi_a)V\|_2= \sqrt{\frac{1-\ln 2}{a\pi}}\|v\|_2$.
\item $\T(\psi_a)V(W^{1,2}(\RR)) \subseteq W^{1,2}(\RR)$.
\item $V\T(\psi_a)V \in \mathfrak S_2(H)$ and $\|V\T(\psi_a)V\|_2\le \frac{\pi+1}{\pi\sqrt2}\|v\|_2^2$.
\item Given $\eps> 0$, there is $\l_\eps\in\CC\setminus\RR$ such that $\|V(\l_\eps - A)^{-1}\| <\eps$.
\een
\el

\bpf
Observe that for any $h\in L^1\cap L^2$ we have
\begeq\label{hv}
(\T(h)V)x(s) = \int_\RR h(t)v(s-t)x(-s+2t)dt.
\eq
Hence, $\T(h)V \in \mathfrak S_2(H)$ and 
$$\|\T(h)V\|_2^2 = \frac12\int_\RR \int_\RR |h(t)v(s-t)|^2 ds dt = \frac12\|h\|_2^2\|v\|_2^2.
$$
Plugging in the norms $\|\varphi_a\|_2$ and $\|\psi_a\|_2$ from Examples \ref{trapez} and \ref{ex2} establishes Properties 1 and 2. 

To prove Property 3, pick $z > 0$ and let $R = R(z; A) = (z-A)^{-1}$. Using the 
definition of the generator of a Banach module, we have
\begeq\label{res}
Rx = \int_\RR f_z(t)T(-t)xdt, x\in H, 
\eq 
where $\widehat f_z(\l) = (\l-z)^{-1}$.
Then for any $h\in L^1\cap L^2$, letting $h_t = T(t)h$, one easily gets
\[
(\T(h)V)Rx = R(\T(h_t)V)x,
\]
after plugging in \eqref{qq} and \eqref{res}. Hence, Property 3 follows.

Next, observe that for any $h\in L^1\cap L^2$ we have
\begeq\label{hv1}
V(\T(h)V)x(s) = \int_\RR v(s)h(t)v(-s-t)x(s+2t)dt.
\eq
Hence, $V\T(h)V \in \mathfrak S_2(H)$ and 
$$\|V\T(h)V\|_2^2 = \frac12\int_\RR \int_\RR |v(s)h(t)v(-s-t)|^2 ds dt \le \frac12\|h\|_\infty^2\|v\|_2^4,
$$
and the estimate for   $\|\psi_a\|_\infty$ from Example  \ref{ex2} yields Property 4. 

Finally, observe that \eqref{res} yields
\[
VR(\l_\eps; A)x(s) = \int_\RR v(s)f_{i\l_\eps}(t)x(-s-t)dt,\ x\in H, \l_\eps\in\RR\setminus\{0\},
\]
and, hence,  
\[
\|VR(\l_\eps; A)\|_2^2 = \frac1{2\pi}\|v\|_2^2\int_\RR\frac{dt}{|t-i\l_\eps|^2} = \frac1{2\l_\eps}\|v\|_2^2,
\]
which implies Property 5.
\epf

From the above lemma, it is clear that we can choose $a > 0$ such that $\|\T(\psi_a)V\|_2< 1$. Then operator $U = I+\T(\psi_a)V \in B(H)$ is invertible and the estimates in the lemma together with \eqref{eqmet2} allow us to use
\cite[Theorem 3.3]{BKU19} to obtain the following result.

\bt\label{simt}
Consider an operator $\mathscr L$ with $V$ of the form \eqref{opv} and the functions $\varphi_a$ and $\psi_a$, $a >0$,  defined in Examples \ref{trapez} and \ref{ex2}. 
Pick $a = 4\frac{1-\ln 2}\pi\|v\|_2^2$. Then
$\|\T(\psi_a)V\|_2= \frac12$,  $U = I+\T(\psi_a)V \in B(H)$, $U^{-1}\in B(H)$, and $\|U^{-1}-I\|_2 \le 1$. Moreover, $U$ is the similarity transform of $\mathscr L$ into $-i\frac d{dt} - B$, where
\[
\bs
B & = \T(\varphi_a)V + U^{-1}(V\T(\psi_a)V - (\T(\psi_a)V)\T(\varphi_a)V),\\
& = U^{-1}(V\T(\psi_a)V + \T(\varphi_a)V) \in \mathfrak S_2(H),
\end{split}
\]
and we have $\|B\|_2 \le \frac{\sqrt2}\pi\left(4{\sqrt{\frac{1-\ln 2}3}}+\pi+1\right)\|v\|_2^2\le 2.45\|v\|_2^2$.
\et

\bpf
Even though the assumptions of \cite[Theorem 3.3]{BKU19} are slightly different, its proof applies nearly verbatim to establish the similarity of  $\mathscr L$ and $-i\frac d{dt} - B$. The postulated estimates are then easily obtained by direct computation.
\epf

Theorems \ref{hspert} and \ref{simt} immediately yield 
Theorem \ref{maint}.
%
%

\brem
We observe that analogs of Theorem \ref{maint} hold for any self-adjoint operator $A$ and a perturbation $V\in \mathfrak L_A(H)$ for which the properties of Lemma \ref{ml} hold without the specific estimates of the Hilbert-Schmidt norms. Moreover, Properties 1, 2, and 4 may be replaced by the following weaker assumptions: 
\begin{itemize}
\item $\T(\psi_a)V \in B(H)$ and there is $a > 0$ such that $\|\T(\psi_a)V\| < 1$.
\item $V\T(\psi_a)V + \T(\varphi_a)V \in \mathfrak S_2(H)$.
\end{itemize}
\erem

\section{Appendix.}\label{app}

In this section, we collect the proofs that we include for completeness of the exposition.

\bpf[Proof of Theorem \ref{hspert}]
Let $E_n = E([-n,n])$ be the spectral projection corresponding to $A$ and the interval $[-n, n]$, $n\in\NN$, and $\tE_n = I - E_n$. Similarly, let $A_n = E_n A = AE_n$, $\tA = A - A_n$, $B_n = E_nBE_n$ and $\tB_n = B-B_n$. Observe that $\|B\|_2^2 = \|B_n\|_2^2+\|\tB_n\|_2^2$, and the sequence $(b_n)$ with $b_n = \|\tB_n\|_2$, $n\in\NN$, is in $\ell^2(\NN)$. Observe also that for $\l\in\CC\setminus\RR$ such that $|\Re e\l | \ge n +2\|B\|_2$ or  $|\Im m\l | \ge 2\|B\|_2$ we have
\[
\begin{split}
(\l - A -B_n)^{-1} & = 
(\l - A)^{-1}\tE_n + E_n(\l - A -B_n)^{-1}E_n \\
& = (\l - A)^{-1}\left(\tE_n + \sum_{k = 0}^\infty \left(B_n(\l - A_n)^{-1}\right)^k\right),
\end{split}
\] 
where the series converges absolutely since $\left\|\left(B_n(\l - A_n)^{-1}\right)\right\| \le \frac12$ due to ${\mathrm{dist } (\l, \s(A_n)) > 2\|B\|_2 \ge 2\|B_n\|}$.
Hence for $\l \in \CC\setminus \RR$ with ${\mathrm{dist }} (\l, [-n,n]) > 2\|B\|_2$ we have
\[
\left\|(\l - A -B_n)^{-1}\right\|  \le
\frac1{|\Im m \l |} \left(1 + \sum_{k = 0}^\infty 2^{-k}\right) = \frac3{|\Im m \l |}.
\] 
For any $n\in\NN$ let 
\[
\begin{split}
Q_n = \big\{\l\in \CC:\ &{|\Im m \l |} > 3\|\tB_n\|_2 
\mbox{ and } \\
&{|\Im m \l |} > 2\|B\|_2,  \mbox{ if } |\Re e\l | \le n+ 2\|B\|_2\big\}.
\end{split}
\]
Then for any $\l\in Q_n$  we have
\[
\begin{split}
(\l - A -B)^{-1} & = (\l - A-B_n)^{-1}\left(I - \tB_n(\l - A-B_n)^{-1}\right)^{-1} \\
& = (\l - A -B_n)^{-1}\left(\sum_{k = 0}^\infty \left(\tB_n(\l - A-B_n)^{-1}\right)^k\right)\in B(H),
\end{split}
\] 
and the result follows by considering the union of $Q_n$, $n\in\NN$.
\epf

\brem
We note that for an explicitly known operator $B$ the above proof essentially yields an algorithm for constructing a function $f\in L^2$ that envelops the spectrum $\s(A+B)$.
\erem

\bpf[Proof of Lemma \ref{funest}]
Observe that for any $a > 0$ we have
\[
\bs
\|f\|_1 & =\int\left|\frac1{a+it}(a+it)f(t)\right|dt
\le \left(\int\frac{dt}{a^2+t^2}\right)^{\frac12} \cdot\left(\int|(a+it)f(t)|^2dt\right)^{\frac12}  \\
& \le \sqrt{\frac\pi a}\left(a\|f\|_2 +\frac1{\sqrt{2\pi}}\|\widehat f'\|_2\right)
= \frac1{\sqrt 2}\left(\sqrt a\|\widehat f\|_2 +\frac1{\sqrt a}\|\widehat f'\|_2\right),
\end{split}
\]
by the Cauchy-Schwarz inequality. Plugging in $a = \frac{\|\widehat f'\|_2}{\|\widehat f\|_2}$, yields the desired result.
\epf

{\subsection* {Acknowledgement}}

The first and  third authors were supported in part by the RFBR grant 19-01-00732.

\bibliographystyle{siam}
\bibliography{../refs}

\end{document}